\input amstex
\documentstyle{amsppt}
\magnification 1100
\NoRunningHeads
\NoBlackBoxes
\document

\def\Q{\Bbb Q}

\def\h{\frak h}

\def\ell{{\text{ell}}}

\def\m{\frak m}
\def\O{\Cal O}

\def\R{\Bbb R}
\def\h1{\hat{\bold 1}}

\def\Hom{\text{Hom}}

\def\Ua{U_q(\tilde\g)}
\def\U2{{\Ua}_2}
\def\g{\frak g}

\def\Z{\Bbb Z}

\def\<{\langle}
\def\>{\rangle}
\def\o{\otimes}
\def\e{\varepsilon}

\def\End{\text{End}}

\def\bA{\overline{A}}
\topmatter
\title On Finite-Dimensional Semisimple and Cosemisimple 
Hopf Algebras in Positive Characteristic\endtitle
\author {\rm Pavel Etingof and Shlomo Gelaki \linebreak
Department of Mathematics\linebreak
Harvard University\linebreak 
Cambridge, MA 02138, USA\linebreak
e-mail: etingof\@math.harvard.edu\linebreak shlomi\@math.harvard.edu}
\endauthor
\endtopmatter
\centerline{July 8, 1998}
\head {\bf Introduction}\endhead
Recently, important progress has been made in the study of
finite-dimensional semisimple Hopf algebras over a field of characteristic
zero (see [Mo] and references therein). Yet, very little is
known over a field $k$ of positive characteristic.
In this paper we first prove in Theorem 2.1 that any finite-dimensional 
semisimple and cosemisimple Hopf algebra over $k$ can be lifted to a 
semisimple Hopf algebra of the same dimension over a field of 
characteristic zero. Moreover, we prove in Theorems 2.2 and 2.3 that this 
lifting is functorial, 
and that it carries a (quasi)triangular object to a (quasi)triangular object.
We then use these lifting theorems to
prove some results on finite-dimensional semisimple and 
cosemisimple Hopf algebras $A$ over $k,$
notably Kaplansky's 5th conjecture from 1975 on the order of the antipode 
of $A$ \cite{K}. These results have already been proved over a
field of characteristic zero,
so in a sense we demonstrate that it is sufficient to consider semisimple
Hopf algebras over such a field (they are also cosemisimple
\cite{LR2}), and then to use our Lifting Theorems 2.1, 2.2 and 2.3 to prove 
them for semisimple and cosemisimple Hopf algebras over a field of positive
characteristic. In our proof of Lifting Theorems 2.1 and 2.2 we use standard
arguments of deformation theory from positive to zero characteristic. The
key ingredient of the proof is the theorem that the bialgebra cohomology
groups vanish. We conclude the paper by proving in Theorem 4.2 that any 
semisimple Hopf algebra $A$ of dimension $d>2$ over a field $k$ of 
characteristic $p>d^{\varphi(d)/2}$ (here $\varphi$ is the Euler function), 
is also cosemisimple. This result was known in characteristic $0$ 
\cite{LR2}, and in characteristic $p>(2d^2)^{2d^2-4}$ \cite{So}.

\head {\bf 1. The bialgebra cohomology}\endhead
In the proof of our lifting theorems we will use the bialgebra cohomology
\cite{GS}. Let $A,B$ be
bialgebras over any field $F,$ and let $\phi:A\to B$ be a 
homomorphism of bialgebras. We define 
$C^{p,q}=C^{p,q}(A,B,\phi)=\Hom_F(A^{\o(p+1)},B^{\o(q+1)})$, 
$p,q\ge 0,$ and two differentials, the algebra differential
$d_a^{p,q}:C^{p,q}\to C^{p+1,q}$ and the 
coalgebra differential 
$d_c^{p,q}: C^{p,q}\to C^{p,q+1},$ by the following formulas:
$$
\gather
(d_a^{p,q}f)(a_1\o \cdots \o a_{p+2})=\\
(-1)^{p+1}\phi^{\o (q+1)}(\Delta_{q+1}(a_1))f(a_2\o \cdots \o a_{p+2})
+(-1)^{p+2}f(a_1a_2\o \cdots \o a_{p+2})+\cdots \\ +f(a_1\o 
\cdots \o a_{p+1}a_{p+2})-
f(a_1\o \cdots \o a_{p+1})\phi^{\o (q+1)}(\Delta_{q+1}(a_{p+2})),\tag 1.1
\endgather $$
$$
\gather
(d_c^{p,q}f)(a_1\o \cdots \o a_{p+1})=\\
\sum \phi(a_1^{(1)}\cdots a_{p+1}^{(1)})\o f(a_1^{(2)}\o \cdots \o 
a_{p+1}^{(2)})-\\
-(\Delta\o I^{\o q})(f(a_1\o \cdots \o a_{p+1}))+\cdots 
+(-1)^{q+1} (I^{\o q}\o \Delta)
(f(a_1\o \cdots \o a_{p+1}))+\\
+(-1)^{q+2}\sum f(a_1^{(1)}\o 
\cdots \o a_{p+1}^{(1)}) \o \phi(a_1^{(2)}\cdots a_{p+1}^{(2)}),\tag 1.2
\endgather
$$
where $\Delta(a)=\sum a^{(1)}\o a^{(2)},$ $\Delta_q:A\to A^{\o q},$ 
$q\ge 2,$ is the iterated coproduct ($\Delta_1=I$ is the identity map).
It is straightforward to check that $d_a^2=0,d_c^2=0$ and $d_ad_c=d_cd_a$. 
Thus, $(C^{\bullet,\bullet},d_a,d_c)$ forms a bicomplex. Consider the 
corresponding 
total complex with $C^n(A,B,\phi)=\oplus_{p+q=n}C^{p,q}$, and
the differential $d$ determined by
$d_{|C^{p,q}}=d_a^{p,q}+(-1)^pd_c^{p,q}$. 
We call the total complex $(C^{\bullet}(A,B,\phi),d)$ the bialgebra cochain 
complex of $(A,B,\phi).$ 
The cohomology of $C^{\bullet}(A,B,\phi)$ is called the bialgebra 
cohomology of 
$(A,B,\phi)$ and denoted by $H^{\bullet}(A,B,\phi)$. In the special case in 
which $A=B$ and $\phi=I,$ we write 
$C^{\bullet}(A),H^{\bullet}(A)$ etc. 

\proclaim{Theorem 1.1} {\bf [St]} If $A$ is a finite-dimensional 
semisimple and cosemisimple Hopf algebra over any field, 
then $H^{\bullet}(A)=0.$
\endproclaim

This theorem has the following generalization.

\proclaim{Theorem 1.2} Let $A$ be a finite-dimensional semisimple 
Hopf algebra and $B$ a finite-dimensional cosemisimple Hopf algebra over any 
field. If $\phi:A\to B$ is a homomorphism of Hopf algebras, then 
$H^{\bullet}(A,B,\phi)=0.$ \endproclaim

\demo{Proof} 
Enlarge the bicomplex $C^{p,q}$ to ${\widehat C}^{p,q}$ by setting 
${\widehat C}^{p,q}=C^{p,q}$ for $p,q\ge 0,$ ${\widehat C}^{-1,q}=B^{\o 
(q+1)}$ and ${\widehat C}^{p,q}=0$ otherwise, and taking the 
differentials $d_a,d_c$ determined 
in (1.2). Let ${\widehat C}^{\bullet}={\widehat C}^{\bullet}(A,B,\phi)$ 
denote the corresponding total complex. Also set $D^q(B)=C^{-1,q}(B)$ and 
let $(D^{\bullet}(B), d_c)$ be the corresponding coalgebra complex. One has 
the following exact sequence of complexes: 
$$ 
0\to C^n(A,B,\phi)\to {\widehat C}^n(A,B,\phi)\to D^{n+1}(B)\to 0,
$$
and thus the following long exact sequence of cohomology:
$$
\cdots\to H^{i+1}(D^\bullet(B))\to H^i(C^\bullet(A,B,\phi))\to H^i({\widehat 
C}^{\bullet}(A,B,\phi))\to H^{i+2}(D^\bullet(B))\to \cdots.
$$

A theorem of Hochschild [Ho, Theorems 3.1 and 4.1] states 
that the cohomology of a finite-dimensional semisimple algebra 
over an algebraically closed field with coefficients in any bimodule 
vanishes in positive dimensions. In particular, this implies 
that $H^i(D^\bullet(B))=0$ for all $i,$ 
and hence that $H^{\bullet}(C^\bullet (A,B,\phi))=
H^{\bullet}({\widehat C}^\bullet (A,B,\phi)).$ Thus it suffices to show that 
$H^{\bullet}({\widehat C}^\bullet (A,B,\phi))=0.$
 
We shall apply the following standard lemma in homological algebra (a dual 
version of [Lo, Lemma 1.0.12]):
Let $E^{\bullet,\bullet}$ be a cochain double complex with $E^{p,q}=0$ 
unless $q\ge 0,$ $p\ge -1,$ and let $E^{\bullet}$ be the corresponding total 
complex. Suppose that the columns of $E^{\bullet,\bullet}$ are acyclic 
except in degree $p=-1.$ Then $H^i(E^{\bullet})=H^i(K^{\bullet}),$ where 
$K^q=ker(d:E^{-1,q}\to E^{0,q}).$

Apply the above to $E^{p,q}={\widehat C}^{p,q}.$ In this case $K^q=(B^{\o 
(q+1)})^A$ (the space of elements of $B^{\o (q+1)}$ which commute with 
$\Delta_{q+1}(\phi(a))$ for all $a\in A$). 
Indeed, the operator $d:E^{-1,q}=B^{\o(q+1)}\to E^{0,q}=A^*\o B^{\o(q+1)}$ is 
given by $d(b_1\o \cdots \o b_{q+1})(a)=[\phi^{\o 
(q+1)}(\Delta_{q+1}(a)),b_1\o 
\cdots \o b_{q+1}].$ Since the $q$-th column of $E^{p,q}$ is the Hochschild 
complex of $A$ with coefficients in $B^{\o(q+1)}$ (with degree shifted 
down by $1$), by Hochschild's theorem it is acyclic in $p\ge 0.$ 
Therefore, $H^i({\widehat C}^{\bullet}(A,B,\phi))=H^{i}((B^{\bullet 
+1})^A)$ which implies that
$$
\gather
H^i(C^{\bullet}(A,B,\phi))=H^{i}(D^{\bullet}(B)^A).
\tag 1.3
\endgather $$

Let $H=Im(\phi),$ and $\phi_i:H\to B$ be the 
corresponding injective Hopf algebra map. Notice that $H$ 
is both semisimple and cosemisimple. Now, since $H=Im(\phi)$ it follows by 
(1.3) that $H^i(C^{\bullet}(A,B,\phi))=H^{i}(D^{\bullet}(B)^H).$ By 
(1.3), this equals $H^i(C^{\bullet}(H,B,\phi_i)),$ which equals by duality 
to $H^i(C^{\bullet}(B^*,H^*,\phi_i^*)).$ Finally, by (1.3), this equals 
to $H^{i}(D^{\bullet}(H^*)^{B^*}),$ and since $H^*=Im(B^*),$ this
equals $H^{i}(D^{\bullet}(H^*)^{H^*})$ which equals $0$ by Theorem 1.1.
$\square$\enddemo

As an example, let us prove a simple corollary of this theorem. 

\proclaim{Corollary 1.3} Let $A$ be a finite-dimensional semisimple Hopf 
algebra and $B$ a finite-dimensional cosemisimple Hopf algebra over any 
field. Then $Hom_{Hopf}(A,B)$ is finite. 
In particular the group of Hopf automorphisms of a finite-dimensional 
semisimple and cosemisimple Hopf algebra $A$ is finite. 
\endproclaim

\demo{Proof} 
Consider the variety $Hom_{Hopf}(A,B).$ The Zariski tangent space to a point 
$\phi\in Hom_{Hopf}(A,B)$ equals 
$H^0(A,B,\phi)$ which equals $0$ by Theorem 1.2. This implies that 
$Hom_{Hopf}(A,B)$ is finite (see [Sh, Ex. 3, p.54 and Th. 3, p.78]), 
and the corollary is proved.
(Notice that to prove the finiteness of the automorphism group 
it is enough to use Theorem 1.1, which states that the Lie algebra of this 
group is zero; so in this case Theorem 1.2 is not needed). 
$\square$\enddemo

\proclaim{Remark 1.4} {\rm Let us give an alternative proof to the fact 
that $Aut_{Hopf}(A)$ is finite, which does not use cohomology. Let $f$ be 
a bialgebra derivation of the Hopf algebra $A.$
Since $A$ is semisimple $f$ is inner. Let $x\in A$ be such that $f(a)=xa-ax$
for all $a\in A.$ Note that $\Delta(f(a))=(I\o f+f\o I)\Delta(a)$ for all
$a\in A,$ is equivalent to $[\Delta(x)-1\o x-x\o 1,\Delta(a)]=0$ for all
$a\in A.$ By Drinfeld's theorem [D], if $\lambda$ is an integral of $A^*$
and $b\in A\o A$ is such that $[b,\Delta(a)]=0$ for all $a\in A$, then 
$(\lambda\o I)(b)$ is central in $A$. Applying this to 
$b=\Delta(x)-x\o 1-1\o x$ we get that $(\lambda\o I)(\Delta(x)-1\o x-x\o 1)
=-\lambda(1)x$ is central in $A$. Since $A$ is cosemisimple
we can choose $\lambda$ such that 
$\lambda(1)\ne 0$. Hence $x$ is central and $f=0.$ 
This implies that the 
Lie algebra of the group of Hopf automorphisms is zero, and hence that it  
is finite. $\square$}
\endproclaim

\proclaim{Corollary 1.5} Let $A$ be a finite-dimensional cosemisimple 
Hopf algebra over any field. Then $A$ admits finitely many  
quasitriangular structures. 
\endproclaim

\demo{Proof} A quasitriangular structure on $A$ determines a 
homomorphism of Hopf algebras $A^{*cop}\to A,$ where $A^{*cop}$ is the Hopf 
algebra obtained from $A^*$ by taking the opposite comultiplication. Therefore 
the result follows from Corollary 1.3. 
$\square$\enddemo

Corollaries 1.3,1.5 were proved by Radford [R1,R2] in characteristic $0.$
He also proved them in characteristic $p$ bigger than the dimension of the 
Hopf algebra, under some additional assumptions. Schneider gave a proof 
of the finiteness of $Aut_{Hopf}(A)$ which is the same as our proof 
\cite{Sc}, and Waterhouse proved the same result without using cohomology 
\cite{W}.

\head {\bf 2. The Lifting Theorems}\endhead
Let $p$ be a prime number and $k$ be a perfect field of 
characteristic $p.$ Let $\O=W(k)$ be the ring of Witt vectors of $k$ [Se,
Sections 2.5, 2.6], and $K$ the field of fractions of $\O.$ 
Recall that $\O$ is a local complete discrete valuation ring, and
that the characteristic of $K$ is zero. Let $\m$
be the maximal ideal in $\O,$ which is generated by $p.$ One has 
$\m^n/\m^{n+1}=k$ for any $n\ge 0$ (here $\m^0=\O$). In the case $k=F_p$
one has that $\O=\Z_{p}$ is the ring of $p-$adic integers, and
$K=\Q_{p}.$

\proclaim{Theorem 2.1} Let $A$ be a semisimple and cosemisimple 
$N$-dimensional Hopf algebra over an algebraically closed field $k$ of 
characteristic $p.$ Then:

(i) There exists a unique (up to isomorphism) 
Hopf algebra $\bA$ over $\O$ which is free of rank $N$ as an 
$\O$-module, and such that $\bA/p\bA$ is isomorphic to $A$
as a Hopf algebra. 

(ii) The $K$-Hopf algebra $A_0=\bA\o_\O K$
is semisimple and cosemisimple. 
All irreducible $A_0$-modules and comodules over $\overline{K}$ are defined 
already over $K$.
The dimensions 
of irreducible modules and comodules 
and the Grothendieck rings of the categories 
of modules and comodules for $A$ are the same as for $A_0$. 
   
\endproclaim

\demo{Proof} 
(i) To show the existence of $\bA$ we will prove that  
there exists a sequence of
 Hopf algebras $A_n$ over $\O_n:=\O/p^n$ which are free modules
of rank $N$ over $\O_n$ such that $A_1=A$ and $A_{n+1}/p^n$ is 
isomorphic to $A_n$ as a Hopf 
algebra. Then we can fix isomorphisms $f_n: A_{n+1}/p^n\to A_n$ 
and define $\bA$ as $\underleftarrow{\lim}A_n$. 

Our proof is by induction on $n$. The case $n=1$ is clear. Suppose 
we have constructed $A_i$ for $i\le n$, and let us construct 
$A_{n+1}$. Take $A_{n+1}$ to be a free rank $N$ module over $\O_{n+1}$ 
and fix a module isomorphism $f_n:A_{n+1}/p^n\to A_n$. 
Let $E$ be the set of extensions of the product and coproduct of $A_n$
to $A_{n+1}$ (just as module maps). That is, $E$ is the set of pairs 
$(m',\Delta')\in \Hom_{\O_{n+1}}(A_{n+1}^{\o 2}, A_{n+1})\oplus 
\Hom_{\O_{n+1}}(A_{n+1}, A_{n+1}^{\o 2})$ such that $m',\Delta'$ are mapped 
to the product and coproduct of $A_n$ under $f_n$. 

Let $A_{ij}=\Hom_{k }(A^{\o i},A^{\o j})$. We have an action 
of the additive group $A_{21}\oplus A_{12}$ on $E$ which is defined 
as follows. Let $(\mu,\delta)\in A_{21}\oplus A_{12}$. 
Lift them in any way to $(\mu',\delta')\in 
\Hom_{\O_{n+1}}(A_{n+1}^{\o 2}, A_{n+1})\oplus 
\Hom_{\O_{n+1}}(A_{n+1}, A_{n+1}^{\o 2})$. Now define 
$(\mu,\delta)\circ (m',\Delta')=(m'+p^n\mu',\Delta'+p^n\delta')$.
Clearly, this does not depend on the lifting, so it defines a desired 
group action. It is clear that the constructed group action is free. It is 
also easy to
see that it is transitive. Indeed, if $(m',\Delta'),(m'',\Delta'')$ are
two elements of $E$ then $(m'-m'',\Delta'-\Delta'')$ is zero modulo $p^n$,
so it has the form $p^n(\mu',\delta')$.
Let $\mu\in A_{21}$, $\delta\in A_{12}$
be the reductions of $\mu',\delta'$ mod $p$. Then $(m'',\Delta'')=
(\mu,\delta)\circ (m',\Delta')$. Thus, $E$ is a principal homogeneous 
space of $A_{21}\oplus A_{12}$. 

Now define a map $c:E\to A_{31}\oplus A_{22}\oplus A_{13}$, which measures 
the failure of $(m',\Delta')\in E$ to satisfy Hopf algebra axioms. 
This map is defined by the following rule: Consider the element 
$$
a(m',\Delta')= (m'(I\o m')-m'(m'\o I),
\Delta'm'-m_{13}'m_{24}'(\Delta'\o\Delta'),
(I\o\Delta')\Delta'-(\Delta'\o I)\Delta').
$$
It is clear that $a$ is zero modulo $p^n$. So there exists $b$ such that
$a=p^nb$. The element $b$ is not unique but unique modulo $p$, so we set 
$c(m',\Delta')=b\text{ mod }p$.  

Observe that $c$ takes values in $C^2(A)$ of the bialgebra cochain complex. 
Moreover, it is straightforward to check that $d\circ c=0$. 

By Theorem 1.1, $H^2(A)=0$. This implies that for any $(m',\Delta')\in E$
there exists $(\mu,\delta)\in C^1(A)=A_{21}\oplus A_{12}$ such
that $c(m',\Delta')=d(\mu,\delta)$. It is easy to check that 
$c(x+p^ny)=c(x)+dy$ for $x\in E$ and $y\in A_{21}\oplus A_{12}$, so if we 
set $m''=m'-p^n\mu,\Delta''=\Delta'-p^n\delta$ then we get 
$c(m'',\Delta'')=0.$ 

It remains to show that $(A_{n+1},m'',\Delta'')$ has the unit, counit, and 
antipode which satisfy the axioms of a Hopf algebra, 
and equal the unit, counit, and antipode of $A_n$ modulo $p^n$. 
For unit and counit it is trivial, since it is well known that 
existence of a unit is preserved under algebra deformations.
To show the existence of the antipode, we have to show that 
there exists $S'':A_{n+1}\to A_{n+1}$ 
such that the antipode equation $m''(S''\o id)\Delta''=i\e$ holds, where 
$i$ is the unit and $\e$ the counit. Consider the map 
$T:\End_{\O_{n+1}}(A_{n+1})\to \End_{\O_{n+1}}(A_{n+1})$ given by 
$T(S'')=m''(S''\o id)\Delta''$. This map is a linear isomorphism modulo 
$p$ (since the antipode in a Hopf algebra is unique), so it is a module 
isomorphism. Thus, the antipode equation has a unique solution. 
The uniqueness of this solution implies that it gives the antipode 
of $A_n$ when reduced modulo $p^n$.  
Thus, $A_{n+1}$ is a Hopf algebra which is isomorphic to $A_n$ modulo
$p^n.$ The existence part of the theorem is proved. 

We now prove uniqueness. Let $\bA',\bA''$ be two algebras satisfying the 
conditions of the theorem. Let $A_n',A_n''$ be their reductions modulo
$p^n$. We will show that 
for any isomorphism of Hopf algebras $f_n:A_n'\to A_n''$  
there exists an isomorphism of Hopf algebras $f_{n+1}: A_{n+1}'\to 
A_{n+1}''$ such that $f_{n+1}=f_n$ modulo $p^n$. Since $A_1'=A_1''=A$,
this implies the uniqueness part of the theorem. 

We identify $A_n'$ and $A_n''$ by $f_n$ and assume $A_n'=A_n''$, $f_n=id$.
We also assume that $A_{n+1}',A_{n+1}''$ are the same as modules. 
As we saw, 
the set of extensions of the product and coproduct on $A_n$ one step higher
(with axioms satisfied) is the set of points $(m',\Delta')\in E$ 
such that $c(m',\Delta')=0$. Thus, if we have two such extensions 
$(m',\Delta')$ and $(m'',\Delta'')$ corresponding to $A_{n+1}',A_{n+1}''$ then 
the element $(\mu,\delta)=(m'-m'',\Delta'-\Delta'')\in C^1(A)$ (here we 
use the transitivity of the action of $A_{12}\oplus A_{21}$ on $E$) 
satisfies $d(\mu,\delta)=0$. But by Theorem 1.1, $H^1(A)=0$. Thus, 
$(\mu,\delta)=d\gamma$, where $\gamma\in \End_k(A)$. Let 
$\tilde\gamma$ be an extension of $\gamma$ to a module map 
$A_{n+1}'\to A_{n+1}''$ (i.e. $\tilde\gamma$ equals $\gamma$ modulo $p$), 
and let $\eta=id-p^n\tilde\gamma: A_{n+1}'\to A_{n+1}''$ (it is clear 
that this depends only on $\gamma$ and not on $\tilde\gamma$). 
This is the required isomorphism.

(ii) It remains to show that $A_0$ is semisimple and cosemisimple, and 
that dimensions of irreducible modules and comodules and Grothendieck rings 
are the same.  To do this, we let 
$n_i$ be the dimensions of irreducible representations of $A$.  
Then $A=\oplus M_{n_i}(k )$, where $M_n$ is the 
matrix algebra of size $n$. 
By Hochschild's theorem, the algebra $A$ has zero Hochschild cohomology, 
so it has a unique lifting to $\O$, namely $\oplus M_{n_i}(\O)$. 
Thus, $\bA=\oplus M_{n_i}(\O)$ as an $\O$-algebra. This implies 
that $A_0=\oplus M_{n_i}(K)$ and hence it is semisimple. 
The fact that $A_0$ is cosemisimple is shown in the same way
by taking the duals (it also follows from [LR2]). The 
fact that all $A_0$ modules and comodules over $\overline{K}$ are 
defined over $K$ and the 
identity of dimensions of the irreducible modules
and comodules and of Grothendieck rings is now clear. $\square$\enddemo

In the following we consider lifting of homomorphisms.

\proclaim{Theorem 2.2} Let $A,B$ be finite-dimensional semisimple and 
cosemisimple Hopf algebras over an algebraically closed field $k$ of
characteristic $p,$ $\phi:A\to B$ be a homomorphism of Hopf 
algebras, and $\bar A,\bar B$ the corresponding liftings to $\O.$ Then
there exists a unique homomorphism of Hopf algebras $\bar \phi:\bar A\to 
\bar B$ such that $\bar \phi=\phi$ modulo $p.$
\endproclaim

\demo{Proof} We first prove existence. It is enough to prove
that there exists a sequence of Hopf algebra 
maps $\phi_n:{\bar A}/p^n\to {\bar B}/p^n,$
$n\ge 1,$ such that $\phi_1=\phi$ and $\phi_n=\phi_{n-1}$ modulo 
$p^{n-1}.$ We 
construct it by induction. For $n=1$ there is nothing to construct as 
$\phi_1=\phi_.$ Suppose $\phi_n$ has been constructed, and choose an 
$\O _{n+1}-$map ${\tilde \phi_{n+1}}:{\bar A}/p^{n+1}\to {\bar 
B}/p^{n+1}$ such that ${\tilde \phi_{n+1}}=\phi_n$ modulo $p^n.$
Write ${\tilde \phi_{n+1}}(xy)-{\tilde 
\phi_{n+1}}(x){\tilde \phi_{n+1}}(y)=p^n{\tilde \psi}(x,y)$ with ${\tilde 
\psi}=\psi$ modulo $p,$ for some $\psi:A\o A\to B.$ Also write 
$({\tilde \phi_{n+1}}\o {\tilde 
\phi_{n+1}})(\Delta(x))-
\Delta 
({\tilde \phi_{n+1}}(x))=p^n{\tilde \eta}(x)$ with ${\tilde \eta}=\eta$ 
modulo $p,$ for some $\eta:A\to B\o B.$ The pair $(\psi,\eta)$ is uniquely 
determined by $\tilde \phi,$ and it is straightforward to check that it 
is a 1-cocycle in the total 
complex $C^n(A,B,\phi).$ Hence by Theorem 1.2, there exists a linear map
$\chi:A\to B$ such that $(\psi,\eta)=d\chi.$ Choose a lifting ${\tilde 
\chi}:{\bar A}/p^{n+1}\to {\bar B}/p^{n+1}$ of $\chi,$ and set 
$\phi_{n+1}={\tilde 
\phi_{n+1}}-p^n{\tilde \chi}.$ Then $\phi_{n+1}$ is a Hopf algebra map.

We now prove uniqueness. Let $\bar\phi,\bar\phi':\bar 
A\to\bar B$ be two liftings of $\phi:A\to B$. 
We prove by induction on $n$ that $\phi_n,\phi_n':
\bar A/p^n\to \bar B/p^n$ 
(the reductions of $\bar\phi,\bar\phi'$ modulo $p^n$)
are equal. Indeed, this is clear for $n=1$; suppose it is known 
for $n=m$. Let $\phi_{m+1}-\phi_{m+1}'=p^m\tilde \chi$, where 
$\tilde\chi=\chi$ modulo $p$, where $\chi: A\to B$ is a linear map. 
It is straightforward to check that $\chi$ is a 0-cocycle in 
$C^\bullet(A,B,\phi)$, which implies by Theorem 1.2 that $\chi=0$ and 
thus $\phi_{m+1}=\phi_{m+1}'$, as desired. This concludes the proof of 
the theorem. $\square$\enddemo

\proclaim{Theorem 2.3} Let $(A,R)$ be a finite-dimensional quasitriangular
semisimple and cosemisimple Hopf algebra over an algebraically closed 
field $k$ of characteristic $p.$ Then there exists a unique $\bar R\in \bar 
A\o \bar A$ such that $(\bar A,\bar R)$ is quasitriangular, and $\bar R=R$ 
modulo $p.$ Furthermore, if $(A,R)$ is triangular then so is $(\bar A,\bar 
R).$ 
\endproclaim

\demo{Proof} Let $R=\sum_i x_i\o y_i,$ and let $D(A)$ denote the 
quantum double of $A.$ Since $(A,R)$ is quasitriangular, the map 
$\phi:(D(A),{\Cal R})\to (A,R)$ determined by the formula $\phi(q\o a)=\sum_i 
q(x_i)y_ia$ is a surjective homomorphism of Hopf algebras such that 
$(\phi\o \phi)(\Cal R)=R$ (here $\Cal R$ is the universal R-matrix of 
$D(A)$). Now consider the Hopf algebra $D(\bar A)$. 
This Hopf algebra is a lifting of $D(A)$, so by 
Theorems 2.1 and 2.2, 
there is a unique isomorphism $D(\bar A)\to \overline{D(A)}$ 
which is the identity modulo $p$. We identify 
$D(\bar A)$ with $\overline{D(A)}$ using this isomorphism. Now, 
by Theorem 2.2, $\phi$ can be lifted to a surjective homomorphism of Hopf 
algebras $\bar \phi:D(\bar A)\to \bar A$. Let $\bar {\Cal R}$ be the 
quasitriangular structure of $D(\bar A),$ and define $\bar R=(\bar \phi\o 
\bar \phi)(\bar {\Cal R}).$ It is straightforward to check that $\bar R$ 
is quasitriangular structure of $\bar A$ which equals $R$ modulo $p$ 
(i.e. is a lifting of $R$). This lifting is unique because any lifting 
$\tilde R$ of $R$ defines a homomorphism of Hopf algebras  
$\theta_{\tilde R}:{\bar A}^{*cop}\to {\bar A},$ which is a lifting of 
$\theta_{R}:A^{*cop}\to A$ corresponding to $R,$ and by Theorem 2.2 the 
lifting of $\theta_{R}$ is unique. 

Now, suppose $(A,R)$ is triangular; that is, $(R^{21})^{-1}=R.$ Then $({\bar 
R}^{21})^{-1}$ and $\bar R$ are both quasitriangular structures on $\bar A$ 
which are liftings of $R.$ By uniqueness of lifting they are equal, and 
$\bar R$ is triangular as well. 
$\square$\enddemo

\proclaim{Corollary 2.4} The assignment $A\to A_0$ 
determines a functor between the categories of finite-dimensional 
semisimple and cosemisimple Hopf algebras over $k$ and finite-dimensional 
semisimple Hopf algebras over $K.$ It also determines a functor between 
the categories of quasitriangular Hopf algebras, and those of triangular 
ones.
\endproclaim

\demo{Proof} The first statement follows from Theorems 2.1 and 2.2. As 
for the second statement, let $A,B$ be finite-dimensional semisimple and 
cosemisimple Hopf algebras over $k.$ Suppose $f:(A,R_A)\to (B,R_B)$ is a 
map of quasitriangular Hopf algebras. Then we have to show that $(\bar 
f\o \bar f)(\bar R_A)=\bar R_B.$ Indeed, set $R'_B=(\bar f\o \bar 
f)(\bar R_A).$ Clearly, the associated map $\theta_{R'_B}:\bar B 
^{*cop}\to \bar B$ is a Hopf algebra map which equals $\theta_{R_B}$ modulo 
$p.$ But, $\theta_{\bar R_B}=\theta_{R_B}$ modulo $p$ too, hence by 
Theorem 2.2 $\bar R_B=R'_B.$
$\square$\enddemo

\head {\bf 3. Applications of the lifting theorems}\endhead

The lifting theorems provide a simple way to prove results about semisimple 
and cosemisimple Hopf algebras in characteristic $p$ which are known 
in characteristic $0$. In this section we give a few applications 
of this sort. We start with Kaplansky's 5th conjecture from 1975.

\proclaim{Theorem 3.1} Let $A$ be a finite-dimensional 
semisimple and cosemisimple Hopf algebra over any field $k.$ Then the
square of the antipode of $A$ is the identity. 
\endproclaim

\demo{Proof} In characteristic $0$, the result is known and due to
Larson and Radford
\cite{LR1}. Suppose $k$ is of characteristic $p$. We can
assume that $k$ is algebraically closed. Using Theorem 2.1(i), we can
construct 
the Hopf algebras $\bA$ and $A_0$. By the characteristic $0$ result, 
the square of the antipode is the identity in $A_0$. Thus it is so in 
$\bA\subset A_0$ and hence in $A=\bA/p\bA$. 
$\square$\enddemo

This result was known in characteristic $p>(dimA)^2$ \cite{LR1}.

\proclaim{Corollary 3.2}  Let $A$ be a finite-dimensional Hopf algebra with 
antipode $S$ over any field $k.$ Then: 

\noindent
(i) $A$ is semisimple and cosemisimple if and only if $S^2=I$ and $dimA\ne 
0$ in $k.$

\noindent
(ii) If $A$ is semisimple and cosemisimple and $k$ is algebraically closed, 
then for any irreducible representation $V$ of $A,$ $dimV\ne 0$ in $k.$ 
\endproclaim

\demo{Proof} (i) Suppose $A$ is semisimple and cosemisimple. Then by Theorem 
3.1, $S^2=I,$ and by [LR1, Theorem 2], $dimA\ne 0$ in $k.$ The "only if" 
part follows from [R4], since $tr(S^2)=tr(I)=dimA\ne 0.$

\noindent
(ii) Follows from [La, Theorem 4.3], which state the same as Corollary 3.2 (ii) 
under the assumption that $A$ is finite-dimensional semisimple and $S^2=I$.
$\square$\enddemo

\proclaim{Corollary 3.3} (i) Let $A$ be a finite-dimensional
semisimple {\bf or} cosemisimple minimal quasitriangular Hopf algebra over 
any field $k.$ Then the square of the antipode of $A$ is the identity. 
If $char(k)=p>0,$ then $p$ does not divide $dimA.$
(ii) Let $A$ be a finite-dimensional semisimple triangular Hopf algebra 
over any field $k.$ Then the Drinfeld element $u$ satisfies $u=S(u)$ 
and $u^2=1.$ 
\endproclaim

\demo{Proof} (i) By [R3, Proposition 14], $A$ is also cosemisimple, and 
hence by Theorem 3.1, $S^2=I.$ The last statement follows now from 
Corollary 3.2.
 
\noindent
(ii) Consider the corresponding minimal Hopf subalgebra $A_R$ of $A.$ 
Clearly, $u\in A_R.$ By (i), $(S_{|A_R})^2=I$ and hence $u$ is 
central in $A_R.$ It is straightforward to check that $tr(u)=tr(S(u))$ in 
any irreducible representation of $A_R$ (see e.g. [EG1]), hence by Corollary 
3.2 (or Theorem 4.3 of \cite{L}), $u=S(u)$. 
But $S(u)=u^{-1}$ since $u$ is a grouplike element. 
Thus $u=u^{-1}$ and $u^2=1$. 
$\square$\enddemo

\proclaim{Theorem 3.4}  Let $A$ be a semisimple and cosemisimple Hopf
algebra of prime dimension $p$ over a field $k.$ Then $A$ is commutative
and cocommutative. 
\endproclaim

\demo{Proof} In characteristic zero, the result is known \cite{Z}. 
For positive characteristic (which has to be different from $p$ by
Corollary 3.2), it suffices to assume that $k$ is 
algebraically closed. In this case $A_0$ is commutative and cocommutative
by the characteristic zero result, thus so are $\bA$ and $A$.    
$\square$\enddemo

\proclaim{Corollary 3.5} Let $A$ be a Hopf algebra of prime dimension $p$ 
over a field $k$ with characteristic $q$ such that $q>p.$ Then $A$ is 
commutative and cocommutative.
\endproclaim

\demo{Proof} We follow the ideas of [Z]. Let $G(A),G(A^*)$ denote the 
groups of grouplike elements of 
$A,A^*$ respectively. By \cite{NZ}, $|G(A)|,|G(A^*)|\in\{1,p\}.$ If either
$|G(A)|=p$ or $|G(A^*)|=p$ then $A=k\Z_p$ and we are done. Suppose 
$|G(A)|=|G(A^*)|=1.$ Then $A$ and $A^*$ are unimodular and 
hence by Radford's formula [R5], $S^4=I.$ Suppose that $A$ is 
not semisimple. Then by [R4], $tr(S^2)=0$ in 
$k.$ But since $S^4=I,$ this implies that there exist integers 
$0\le a,b\le p$ such that $a-b=0$ in $k,$ and $a+b=p.$ Since $q>p$ this is 
impossible, and hence $A$ is semisimple. Similarly, $A^*$ is semisimple and 
the result follows from Theorem 3.4. 
$\square$\enddemo

\proclaim{Theorem 3.6}  Let $A$ be a semisimple and cosemisimple Hopf
algebra of dimension $pq$ over a field $k,$ where $p,q$ 
are distinct primes. Then $A$ is commutative or cocommutative. 
\endproclaim

\demo{Proof} In characteristic zero, the result is known \cite{EG2}. 
For positive characteristic (which has to be different from $p$ and $q$ by
Corollary 3.2), it suffices to assume that 
$k$ is algebraically closed. In this case $A_0$ is commutative or
cocommutative by the characteristic zero result, 
thus so are $\bA$ and $A$.    
$\square$\enddemo

\proclaim{Theorem 3.7} Let $A$ be a finite-dimensional 
semisimple and cosemisimple Hopf algebra over an algebraically closed 
field, and let $D(A)$ be the quantum double of $A$. Then the dimension of
any irreducible representation of $D(A)$ divides the dimension of $A$. 
\endproclaim

\demo{Proof} In characteristic zero, the result is known \cite{EG1}. 
For positive characteristic, first note that the double of $A$ is
semisimple and cosemisimple, and that taking the double commutes with
lifting (as we mentioned in the proof of Theorem 2.3). 
By the characteristic zero result, the
statement is true for $D(A_0)\o_K \overline{K}.$ By Theorem 2.1(ii), all
representations of $D(A_0)\o_K\overline{K}$ split already over $K$. Therefore,
applying Theorem 2.1(ii) again we get that 
the result also holds for $D(A)$.    
$\square$\enddemo

\proclaim{Corollary 3.8} Let $A$ be a finite-dimensional 
semisimple and cosemisimple 
quasitriangular Hopf algebra over an algebraically closed field.
Then the dimension of any irreducible representation 
of $A$ divides the dimension of $A$. 
\endproclaim

\demo{Proof} We have a surjective homomorphism of Hopf algebras $D(A)\to
A$ (as we mentioned in the proof of Theorem 2.3), so any irreducible 
$A$-module is also an irreducible $D(A)$-module. Now apply Theorem 3.7. 
$\square$\enddemo

\proclaim{Theorem 3.9} Let $A$ be a finite-dimensional semisimple
and cosemisimple triangular Hopf algebra over any field $k$. 
Then the dimensions of irreducible modules and 
the Grothendieck ring of the category $Rep(A)$ 
coincide with those of some finite group $G$. 
\endproclaim

\demo{Proof} Consider the Hopf algebra $A_0$. Its category of representations 
is equivalent to the one of some finite group $G$ by [EG1, Theorem 2.1].  
But the dimensions of irreducible representations and
the Grothendieck ring are the same for $A$ and $A_0$ by Theorem 2.1. 
$\square$\enddemo

\proclaim{Theorem 3.10} 
Let $A$ be a semisimple and cosemisimple Hopf algebra of dimension $p^n$,
where $p$ is a 
prime, over an algebraically closed field $k$. 
Then $A$ has a non-trivial central grouplike element.
\endproclaim

\demo{Proof} In characteristic zero, the result is known \cite{Ma}. 
For positive 
characteristic (which has to be different from $p$ by Corollary 3.2), 
the statement is true for $A_0\o_K \overline{K}$ 
by the characteristic zero result. 
Let $g$ be a non-trivial central grouplike element of $A_0\o_K
\overline{K}$.
Then $g:A_0^*\to \overline{K}$ is a 1-dimensional representation
such that $(g\o 1)\circ \Delta=(1\o g)\circ\Delta$. But
by Theorem 2.1(ii), $g:A_0^*\to
\overline{K}$
really 
takes values in $K.$ Hence, by Theorem 2.1(ii), this representation is
obtained by lifting of a unique representation $g':A^*\to k$, which 
therefore must satisfy the same 
equation. Then $g'$ is a nontrivial central grouplike element in $A$.   
$\square$\enddemo

\head {\bf 4. When Finite-Dimensional Semisimple Hopf Algebras are 
Cosemisimple}\endhead
 
Let $A$ be any semisimple Hopf algebra of dimension $d>2$ 
over a field $k$ of characteristic $p>d^{\varphi(d)/2}$ (here $\varphi$ is 
the Euler function i.e. $\varphi(d)=\#\{1\le m\le d|(d,m)=1\}$). In this 
section we prove that $A$ is also cosemisimple. 
This result was known in characteristic $0$ \cite{LR2}, and
in characteristic $p>(2d^2)^{2d^2-4}$ \cite{So}. So
our result slightly improves Sommerhauser's result. Our proof is based on 
combining the ideas of Larson and Radford [LR2] with a trivial number 
theoretic lemma.

\proclaim{Lemma 4.1} Let $P=\sum_{m=0}^{r-1} a_mx^m\in\Z[x]$ be a 
polynomial of degree $r-1>1.$ Let $\sum_{m=0}^{r-1}|a_m|=D,$ and suppose 
that $P(e^{2\pi i/r})$ is a non-zero real number. Then for all 
$p>D^{\varphi(r)/2}$ such that $p$ does not divide $r,$ and  
any primitive $rth$ root of unity $\zeta\in \bar F_p$ (the algebraic 
closure of the field of $p$ elements), one has $P(\zeta)\ne 0.$
\endproclaim

\demo{Proof} Let $N=\Pi_{(l,r)=1,l<r/2}P(e^{2\pi il/r})$ be the product 
of all the conjugates of $P(e^{2\pi i/r}).$ It is easy to see that $N\in\Q^*$ 
and is an algebraic integer, so we have that $N$ is a non-zero integer. 
Clearly, $|N|\le D^{\varphi(r)/2},$ and hence $p$ does not divide $N.$ This 
implies that $N$ 
is a unit in the ring $\O=W(\bar F_p)$ of Witt vectors of 
$\bar F_p.$ 
Let $\hat \zeta$ be the lifting of $\zeta$ to 
${\Cal O};$ that is, ${\hat \zeta}^r=1$ and $\hat \zeta=\zeta$ modulo 
$p.$ Then $P(\hat \zeta)$ is a unit in ${\Cal O}$ since 
$\Pi_{(l,r)=1,l<r/2} P({\hat \zeta}^l)=N$ and $P({\hat \zeta}^l)$ is 
an algebraic integer (hence an integer in ${\Cal O}$) for all $l.$ Therefore, 
$P(\zeta)\ne 0$ in $\bar F_p.$ $\square$\enddemo

\proclaim{Theorem 4.2} Let $A$ be a $d-$dimensional Hopf algebra, $d>2,$  
with antipode $S$ over a field $k$ of characteristic 
$p>d^{\varphi(d)/2}.$ Then the following are equivalent:

\noindent
(i) $A$ is semisimple. 

\noindent
(ii) $A$ is cosemisimple.

\noindent
(iii) $S^2=I.$
\endproclaim

\demo{Proof} (i) implies (ii): We may assume that $k$ is 
algebraically closed, and start as in [LR2]. Since $A$ is 
semisimple we have $A=\oplus_i End_k(V_i)$ where $V_i$ are the irreducible 
representations of $A.$ Let ${\Cal I}$ be the set of all $i$'s for which 
$S^2(End_k(V_i))=End_k(V_i)$ (in fact, it is easy to show that all $i$ are 
such, but we do not need it). Since $A$ is unimodular it follows by [R5] 
that for all $a\in A,$ $S^4(a)=gag^{-1}$ for some grouplike element $g\in 
A.$ Thus $S^{4d}=I,$ so if $s$ is the order 
of $S^4$ then $s$ divides $d.$ Let $r=2s$ be the order of $S^2.$ If $r=2,$ 
then $S^4=I$ and the result follows from [LR1] (for $p>d$). If $r=2d$ then 
$A$ is a group algebra, so there is nothing to prove. Therefore we may 
assume $2<r\le d,$ thus $p>r.$ Since for 
all $i\in {\Cal I},$ $S^2_{|End_k(V_i)}$ is an algebra automorphism, it is 
inner; that is $S^2_{|End_k(V_i)}(B)=g_iBg_i^{-1}$ for some $g_i\in 
Aut_k(V_i).$ Since $g_i^r$ is 
central in $End_k(V_i)$ we may assume that $g_i^r=I.$ Let $\zeta$ be a 
primitive $rth$ root of unity in $k,$ and write $V_i=\oplus_{j=0}^{r-1}V_{ij}$ 
where $V_{ij}$ is the eigenspace corresponding to $\zeta^j$ (possibly  
$V_{ij}=0$). Set $d_{ij}=dimV_{ij}.$ Then, 
$$
\gather
T=tr(S^2)=\sum_{i\in {\Cal 
I}}tr(S^2_{|End_k(V_i)})=\sum_{i\in {\Cal I}}tr(g_i)tr(g_i^{-1})=
\sum_{i\in {\Cal I}}(\sum_{j=0}^{r-1}d_{ij}\zeta^j)
(\sum_{j=0}^{r-1}d_{ij}\zeta^{-j}).
\endgather $$ 
By [R4], it is enough to show that $T\ne 0$ in $k.$

Now, since $\zeta^r=1,$ we can write 
$$
\gather
T=\sum_{m=0}^{r-1}a_m\zeta^m,\;\;a_m=\sum_{i\in {\Cal 
I}}(\sum_{j-l=m\,(mod\,r)}d_{ij}d_{il})\in\Z_+.
\endgather $$ 
Consider the polynomial $P(x)=\sum_{m=0}^{r-1}a_mx^m\in 
\Z[x].$ Since in $\Q(e^{2\pi i/r}),$ 
$$
\gather
(\sum_{j=0}^{r-1}d_{mj}e^{2\pi ij/r})
(\sum_{j=0}^{r-1}d_{mj}e^{-2\pi ij/r})=|\sum_{j=0}^{r-1}d_{mj}e^{2\pi 
ij/r}|^2\ge 0,
\endgather $$ 
and $tr(S^2_{|End_k(k)})=1,$ it follows that $P(e^{2\pi i/r})\in 
\R^*.$ Finally, $D=\sum_m a_m\le d$ and 
$\varphi(r)\le \varphi(d)$ (as $r=2d/n,$ where $n\ge 2$ and divides $d$).
Therefore the result follows from Lemma 4.1.

\noindent
(ii) implies (i): Follows from (i) implies (ii) for $A^*.$

\noindent
(i) and (ii) imply (iii): Follows from 
[LR1] where the result is proved for $p>d^2.$ 

\noindent
(iii) implies (i) and (ii): Follows from [R1] since $d\ne 0$ in $k.$ 
$\square$\enddemo

\noindent
{\bf Acknowledgments.} The authors are grateful to David Kazhdan and Murray 
Gerstenhaber for useful discussions and explanations, and to Dima Arinkin 
for reading the manuscript.

\end